\documentclass[11pt]{amsart}
\usepackage{amssymb, latexsym, mathrsfs,  color, mathtools, multirow}
\usepackage[colorlinks=true, pdfstartview=FitV,linkcolor=blue,citecolor=blue,urlcolor=blue]{hyperref}

\newcommand{\GL}{\operatorname{GL}}

\usepackage{fullpage}

\newtheorem {theorem}    {Theorem}[section]

\newtheorem {lemma}      [theorem]    {Lemma}

\theoremstyle{definition}

\numberwithin{equation}{section}

\usepackage[colorinlistoftodos]{todonotes}

\newenvironment{red}{\relax\color{red}}{\hspace*{.5ex}\relax}
\newenvironment{blue}{\relax\color{blue}}{\hspace*{.5ex}\relax}
\newcommand{\ber}{\begin{red}}
	\newcommand{\er}{\end{red}}
\newcommand{\beb}{\begin{blue}}
	\newcommand{\eb}{\end{blue}}

\newcommand{\seteq}{\mathbin{:=}}

\newcommand{\sage}{{\sc SageMath }} 
\newcommand{\mth}{{\sc Mathematica }}
\newcommand{\tf}{{\sc TensorFlow }}
\newcommand{\lgbm}{{\sc LightGBM }}

\begin{document}

\title{Machine-Learning Kronecker Coefficients}

\date{\today}

\author[K.-H. Lee]{Kyu-Hwan Lee$^{\star}$}
\thanks{$^{\star}$This work was partially supported by a grant from the Simons Foundation (\#712100).}
\address{Department of
Mathematics, University of Connecticut, Storrs, CT 06269, U.S.A.}
\email{khlee@math.uconn.edu}

\begin{abstract}
The Kronecker coefficients are the decomposition multiplicities of the tensor product of two irreducible representations of the symmetric group. Unlike the Littlewood--Richardson coefficients, which are the analogues for the general linear group, there is no known combinatorial description of the Kronecker coefficients, and it is an NP-hard problem to decide whether a given Kronecker coefficient is zero or not. 
In this paper, we show that standard machine-learning algorithms such as Nearest Neighbors, Convolutional Neural Networks and Gradient Boosting Decision Trees may be trained to predict whether a given Kronecker coefficient is zero or not. Our results show that a trained machine can efficiently perform this binary classification with high accuracy ($\approx 0.98$).
\end{abstract}

\maketitle


\section{Introduction}

One of the major goals in representation theory is to decompose a representation into irreducible representations, and algebraic combinatorics plays a crucial role as it produces a practical way to describe such a decomposition. A prototype can be found in the decomposition of the tensor product of two irreducible representations of the general linear group $\GL_N(\mathbb C)$, where the Littlewood--Richardson rule completely describes the decomposition using skew semi-standard tableaux. Therefore it comes with surprise to realize that we do not have a similar success with the symmetric group until now.  

To be precise, let $\mathfrak S_n$ be the symmetric group of degree $n$ and consider representations of $\mathfrak S_n$ over $\mathbb C$. The irreducible representations $S_\lambda$ of $\mathfrak S_n$ are parametrized by partitions $\lambda \vdash n$, and  the dimension of $S_\lambda$ is equal to the number of standard tableaux of shape $\lambda$. Consider the tensor product of two irreducible representations $S_\lambda$ and $S_\mu$ for $\lambda, \mu \vdash n$. Then their tensor product is decomposed into a sum of irreducible representations:
\[  S_\lambda \otimes S_\mu = \bigoplus_{\nu\vdash n} g_{\lambda, \mu}^\nu S_\nu  \quad (g_{\lambda, \mu}^\nu \in \mathbb Z_{\ge 0}) . \]
The decomposition multiplicities $g_{\lambda, \mu}^\nu$ are called the {\em Kronecker coefficients}.

In stark contrast with the Littlewood--Richardson coefficients for $\GL_N(\mathbb C)$, no combinatorial description has been known for $g_{\lambda, \mu}^\nu$ since Murnaghan \cite{Mur} initially posed the question in 1938, and it is still considered as one of the main problems in the combinatorial representation theory. For special cases, various combinatorial formulas are known due to Remmel \cite{Rem}, Ballantine--Orellana \cite{BO}, Remmel--Whitehead \cite{RW}, Blasiak--Mulmuley--Sohoni \cite{BMS} and Blasiak \cite{Bla}. Recently, the coefficients $g_{\lambda, \mu}^\nu$ have also been studied from the viewpoint of computational complexity. In particular, it is shown by Ikenmeyer, Mulmuley and Walter \cite{IMW} that deciding whether a given Kronecker coefficient is non-zero is NP-hard, and other interesting results are obtained by Pak and Panova \cite{PP1, PP2}. 

In this paper we apply machine-learning to the Kronecker coefficients $g_{\lambda, \mu}^\nu$ to show that a machine can be trained to learn these coefficients. Though this seems to be the first time for machine-learning to be applied to the coefficients $g_{\lambda, \mu}^\nu$, there has been mounting evidence that a machine can be trained to learn mathematical structures. Indeed, starting with superstring theory \cite{He:2017aed, He:2017set,Krefl:2017yox,Ruehle:2017mzq,Carifio:2017bov}, machine-learning has been applied, for example, to group theory \cite{He:2019nzx}, graph theory \cite{He:2020fdg,Wag}, metric geometry \cite{Ashmore:2019wzb,DLQ}, knot invariants \cite{Jejjala:2019kio}, and to number theory \cite{HLOa, HLOb, HLOc}. Moreover, a new paradigm of AI-aided research in mathematics is emerging \cite{DVB+}.

As the main machine-learning methods for this paper, we utilize Nearest Neighbors, Convolutional Neural Networks and Gradient Boosting Decision Trees (specifically, LightGBM \cite{LGBM}) as they perform better than other standard machine-learning algorithms. For a broad introduction to machine-learning, see \cite{ML,Hastie}. Our datasets of the Kronecker coefficients $g_{\lambda, \mu}^\nu$ are generated by \sage \cite{sage}. 

The results of our experiments clearly show that a machine can be trained to distinguish between nonzero and zero Kronecker coefficients $g_{\lambda, \mu}^\nu$ with high accuracy ($\approx 98\%$). We focus on $\mathfrak S_n$ for $12 \le n \le 14$ since general features are clear enough from these cases, and the sizes of datasets are suitable for our purposes. (When $n=14$, there are $ 2,258,526$ data points in our dataset.) We hope that the results in this paper may invite researchers to further study Kronecker coefficients by taking machine-learning or AI-aided approaches. Specifically, attribute analysis would be the next step, which may produce some conjectures to describe Kronecker coefficients.

This paper is organized as follows. In Section~\ref{s:strategy}, the generation of training data and the experimental set-up are explained. In Section~\ref{s:kron}, we document the experimental outcomes and make some observations. In the final section, we offer some concluding remarks and directions for further research.

\subsection*{Acknowledgments}
We are grateful to Greta Panova for helpful comments.

\section{Methodology}\label{s:strategy}
In this section we explain our experimental set-up. In particular, we present details about how to construct our datasets and overview the machine-learning strategies adopted.

\subsection{Generating Datasets}\label{s:Ltraining}

We keep the notations in the introduction. Additionally, for a partition $\lambda=(\lambda_1 \ge \lambda_2 \ge \cdots) $ of $n$, define $d_\lambda \seteq n- \lambda_1$, called the {\em depth} of $\lambda$. The following theorem suggests a mild restriction for $(\lambda, \mu, \nu)$ to be considered.  

\begin{theorem} \cite[Theorem 2.9.22]{JK} \label{thm-lim}
If $g_{\lambda, \mu}^\nu \neq 0$ then \begin{equation} \label{eqn-dl}  |d_\lambda - d_\mu | \le d_\nu \le d_\lambda + d_\mu . \end{equation}
\end{theorem}
\noindent The inequalities in \eqref{eqn-dl}  provide necessary conditions for $g_{\lambda, \mu}^\nu \neq 0$ and the converse is not true. Actually, when $n=12$, there are 456,533 triples $(\lambda, \mu, \nu)$ in total and 406,919 of them satisfy \eqref{eqn-dl}.  Among them, 280,009 triples give $g_{\lambda, \mu}^\nu \neq 0$ and the remaining 126,910 triples $g_{\lambda, \mu}^\nu=0$. We will consider classification problems for datasets satisfying \eqref{eqn-dl}.  

There are symmetries among $g_{\lambda, \mu}^\nu$ as stated in the lemma below, which will be explicitly taken into account when we make some datasets. 
\begin{lemma} \cite[p.61]{FH} \label{lem-perm}
Let $\lambda, \mu, \nu \vdash n$. Then the Kronecker coefficients $g_{\lambda, \mu}^\nu$ are invariant under the permutations of $\lambda, \mu, \nu$. That is, we have
\[ g_{\lambda,  \mu}^\nu=g_{\mu, \lambda}^\nu=g_{\lambda, \nu}^\mu=g_{\nu, \lambda}^\mu=g_{\mu, \nu}^\lambda=g_{\nu, \mu}^\lambda. \]
\end{lemma}
 
\bigskip

Now, for $n \in \mathbb Z_{>0}$, let $\mathcal P(n)$ be the set of partitions of $n$.
We identify each element $\lambda$ of $\mathcal P(n)$ with a sequence of length $n$ by adding $0$-entries as many as needed. For example, when $n=5$, we have
\begin{align*} \mathcal P(5) = \{ (5,0,0,0,0), (4, 1,0,0,0), (3, 2,0,0,0), (3, 1, 1,0,0), \phantom{0,0), (3, 1, 1,0,0),}\\ (2, 2, 1,0,0), (2, 1, 1, 1,0), (1, 1, 1, 1, 1) \} . \end{align*}

Let $\mathcal Q(n)$ be the set of triples $(\lambda, \mu, \nu) \in \mathcal P(n)^3$ that satisfies \eqref{eqn-dl}. For each triple $(\lambda, \mu, \nu)\in \mathcal Q(n)$, we form the following data points.

(1) a $1d$-array of size $1 \times 3n$ given by concatenation:
\[ v_1(\lambda, \mu, \nu) = [ \lambda , \mu, \nu]. \]

(2) a $2d$-array of size $n \times 3$ given by
\[ v_2( \lambda, \mu, \nu ) = [[\lambda], [\mu], [\nu]],  \] where $[\lambda],[\mu],[\nu]$ are column vectors.

(3) a $3d$-array of size $6 \times n \times 3$ consisting of the permutations of the triple:  
\[ v_3(\lambda, \mu, \nu) = [ [[\lambda], [\mu],[\nu]],[[\lambda],[\nu],[\mu]],[[\mu], [\lambda],[\nu]], [[\mu],[\nu],[\lambda]], [[\nu],[\lambda], [\mu]], [[\nu], [\mu], [\lambda]]] .  \] 
The array $v_3(\lambda, \mu, \nu)$ is defined to reflect Lemma \ref{lem-perm} into the dataset. 

Define a class-valued function $t(\lambda, \mu, \nu)$ on $\mathcal Q(n)$ by
\[ t(\lambda, \mu, \nu) =  \left \{ \begin{array}{ll} 0 & \text{ if $g_{\lambda, \mu}^\nu =0$}, \\ 1 &  \text{ if $g_{\lambda, \mu}^\nu  \neq 0$}. \end{array} \right .     \] 
Finally, we have the following labeled datasets:
\begin{equation}\label{eq.Ldata}
\mathcal{D}_{n}^{(a)}=\{v_a (\lambda, \mu, \nu) \rightarrow t (\lambda, \mu, \nu) \ : \ (\lambda, \mu, \nu ) \in \mathcal Q(n) \} \quad \text{ for $a=1,2,3$.} 
\end{equation}

\subsection{Experimental strategy}\label{sec:strategy}
We perform our experiments following the steps explained below.

\begin{enumerate}

\item Choose $n$. In this paper, we mainly consider $12 \le n \le 14$. When $n=14$, we have $\#(\mathcal Q(14)) = 2,258,526$, which makes the size of dataset reasonably big.

\item Determine an ML classifier. We mostly use Nearest Neighbors (NearN), Convolutional Neural Networks (CNNs) and LightGBM (LGBM).
 
\item Generate datasets of the form $\mathcal{D}_{n}^{(a)}$ as in \eqref{eq.Ldata} using \sage. According to the choice of an ML classifier, we use $a=1$ for NearN and LGBM and $a=2,3$ for CNNs.

\item Form a balanced dataset $\mathcal D$ out of $\mathcal D_{n}^{(a)}$ by sampling the same number of data points from each class $t(\lambda, \mu, \nu)=i$ for $i=0,1$. Choose a subset $\mathcal{T}\subset\mathcal{D}$ and denote its complement by $\mathcal{V}=\mathcal{D}-\mathcal{T}$. We will refer to $\mathcal T$ as the training dataset, and $\mathcal V$ as the validation dataset. We will typically choose 70\% of data points for $\mathcal T$ and 30\% for $\mathcal V$. 
 
\item Train the classifier in step (2) on the set $\mathcal{T}$.  We implement the algorithms using \mth \cite{wolfram}, \tf \cite{TF} and \lgbm \cite{LGBM}.

\item For all triples $(\lambda, \mu, \nu)$ in $\mathcal{V}$, ask the classifier to determine $t(\lambda, \mu, \nu)$, that is, whether  $g_{\lambda, \mu}^\nu =0$ or not. We record the precision and confusion matrix, which together constitute a good measure of performance of the classifier.

\item Repeat steps (4) to (6) for different choices of $\mathcal{T}$. The precision values and confusion matrices recorded below are representative of several repetitions. 

\item Repeat steps (1) to (7) for different choices of $n, a$ and an ML classifier.

\end{enumerate}

\section{Results}\label{s:kron}
In this section we describe our experimental results for the binary classification: $g_{\lambda, \mu}^\nu =0$ or $\neq 0$.

For $n=12,13,14$, we use datasets $\mathcal D^{(1)}_{n}$ to train NearN and LGBM classifiers.  On the other hand, $\mathcal D^{(2)}_{n}$ and $\mathcal D^{(3)}_{n}$ are used for CNN classifiers to utilize the benefit of convolution layers. The results of CNN using $\mathcal D^{(2)}_{n}$ will be denoted by CNN$_2$ and those using $\mathcal D^{(3)}_{n}$ by CNN$_3$. The results are recorded in Tables \ref{t:rank} and \ref{t:conf}. 

\subsection{Specifics}
In the architecture of CNNs, we put one $2d$-convolution layer with 32 filters and one dense layer after flattening the outcome of the convolution layer. Though there are no significant differences in precisions for different kernel sizes of the convolution layers, the following kernel sizes perform slightly better than others: for CNN$_2$, kernel size $(n-4) \times 3$ with input size $n \times 3$, and for CNN$_3$, kernel size $2 \times (n-4)$ with input size $6 \times n \times 3$, where $n=12, 13, 14$. The total numbers of parameters are shown below:
\[ \text{CNN$_2$: $(96 \times n -30)$, \qquad  CNN$_3$: $(192 \times n +866)$.}\]
That is, $1122, 1218, 1314$ for CNN$_2$ and $3170, 3362, 3554$ for CNN$_3$ when $n=12,13,14$, respectively.
The outcomes for CNNs in Table \ref{t:rank} are from these kernel sizes.

The hyper-parameters of LightGBM used in our experiments are below:
{\small \begin{verbatim}             
    'metric': 'auc',            'boosting': 'gbdt',         'num_leaves': 63,       
    'feature_fraction': 0.5,    'bagging_fraction': 0.5,    
    'bagging_freq': 20,         'learning_rate': 0.01
\end{verbatim} }

\subsection{Observations} Let us make some observations from the results recorded in Tables \ref{t:rank} and \ref{t:conf}.
\begin{enumerate}
\item The precisions of different ML classifiers can be compared as
\[ \text{NearN $<$ CNN$_2$ $<$ CNN$_3$ $<$ LGBM}. \]
When $n=14$, LGBM achieved precision 0.9845.   
\item  The classifiers generally perform better as $n$ increases. However, the difference between $n=13$ and $n=14$ is much less, compared to that between $n=12$ and $n=13$. It seems that the precision is being stabilized as $n$ grows. 

\item In the confusion matrix one can see that class 1 ($g_{\lambda, \mu}^\nu \neq 0$) is more often misclassified than class 0 ($g_{\lambda, \mu}^\nu =0$). However, though rarely, the opposite may occur depending on samples. For example, see $n=13$ of CNN$_3$ in Table \ref{t:conf}.
\end{enumerate}

\begin{table}[h!!!]
\begin{center}
{\footnotesize \begin{tabular}{|c|c|c|c|c|c|}
\hline
\multirow{2}{*}{$n$} &\multirow{2}{*}{$\# \mathcal D$} & \multicolumn{4}{|c|}{Precision}\\
\cline{3-6} & & NearN  &CNN$_2$ & CNN$_3$&  LGBM \\ \hline 
$12$ & $ 126,900 \times 2$ & 0.9155 & 0.9529 & 0.9697 & 0.9714  \\
\hline
$13$ & $ 260,000 \times 2$ & 0.9318 & 0.9618 &0.9773& 0.9837  \\
\hline
$14$ & $ 600,000 \times 2$ & 0.9364 & 0.9635 & 0.9772 &  0.9845\\
\hline
\end{tabular}
}
\end{center}
\medskip
\caption{\small \sf
The above table shows the precisions of NearN, CNN$_2$, CNN$_3$ and LGBM classifiers when asked to distinguish $(\lambda, \mu, \nu)$ with $g_{\lambda, \mu}^\nu=0$ from those with $g_{\lambda, \mu}^\nu \neq 0$.
}
\label{t:rank}
\end{table}

\begin{table}[h!!!]
\begin{center}
{\footnotesize \begin{tabular}{|c|c|c|c|c|}
\hline $n$ & NearN & CNN$_2$& CNN$_3$ & LGBM\\ \hline \raisebox{1.85 cm}{$12$}& \includegraphics[width=0.22 \textwidth]{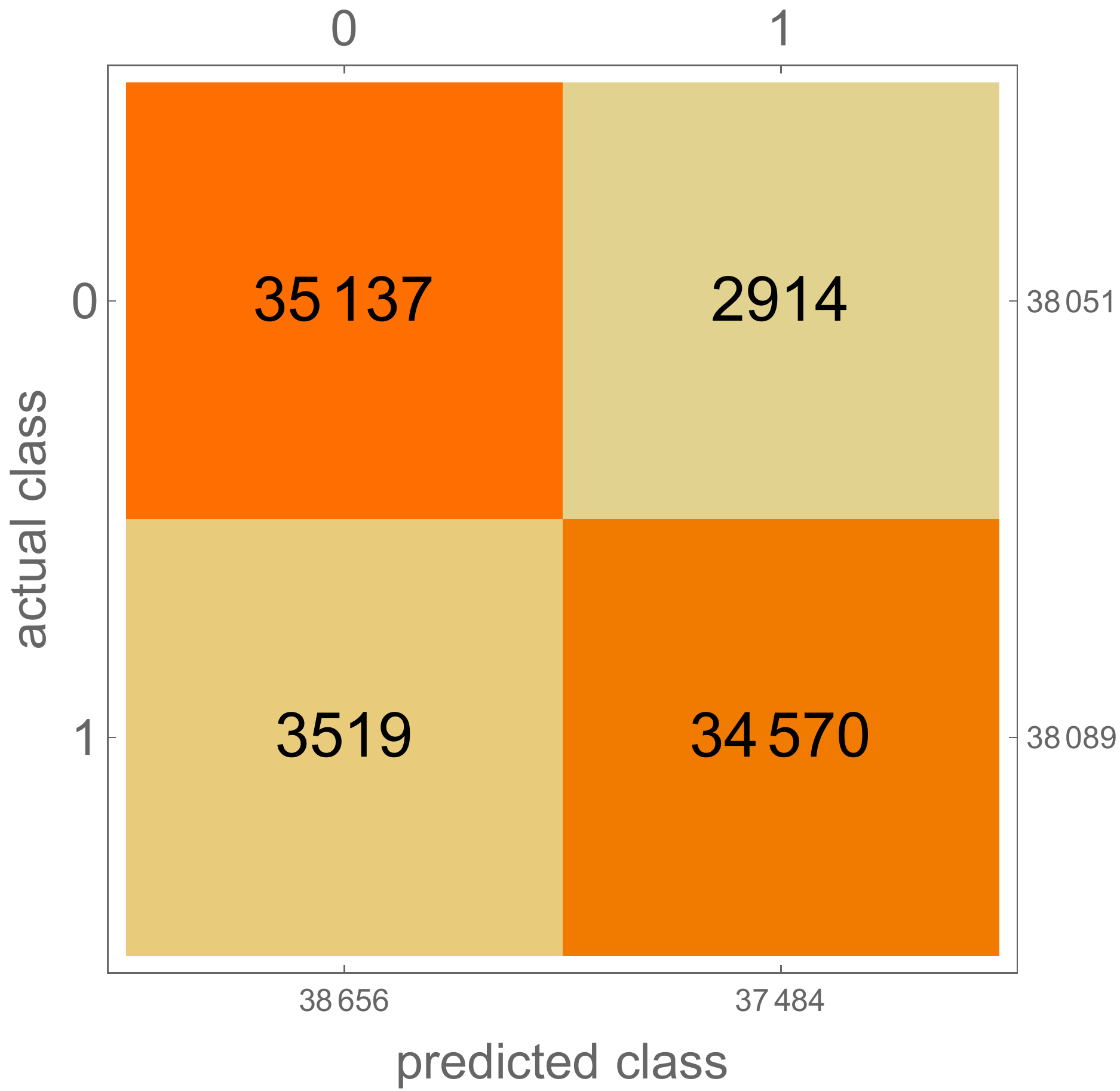}  &\includegraphics[width=0.22 \textwidth]{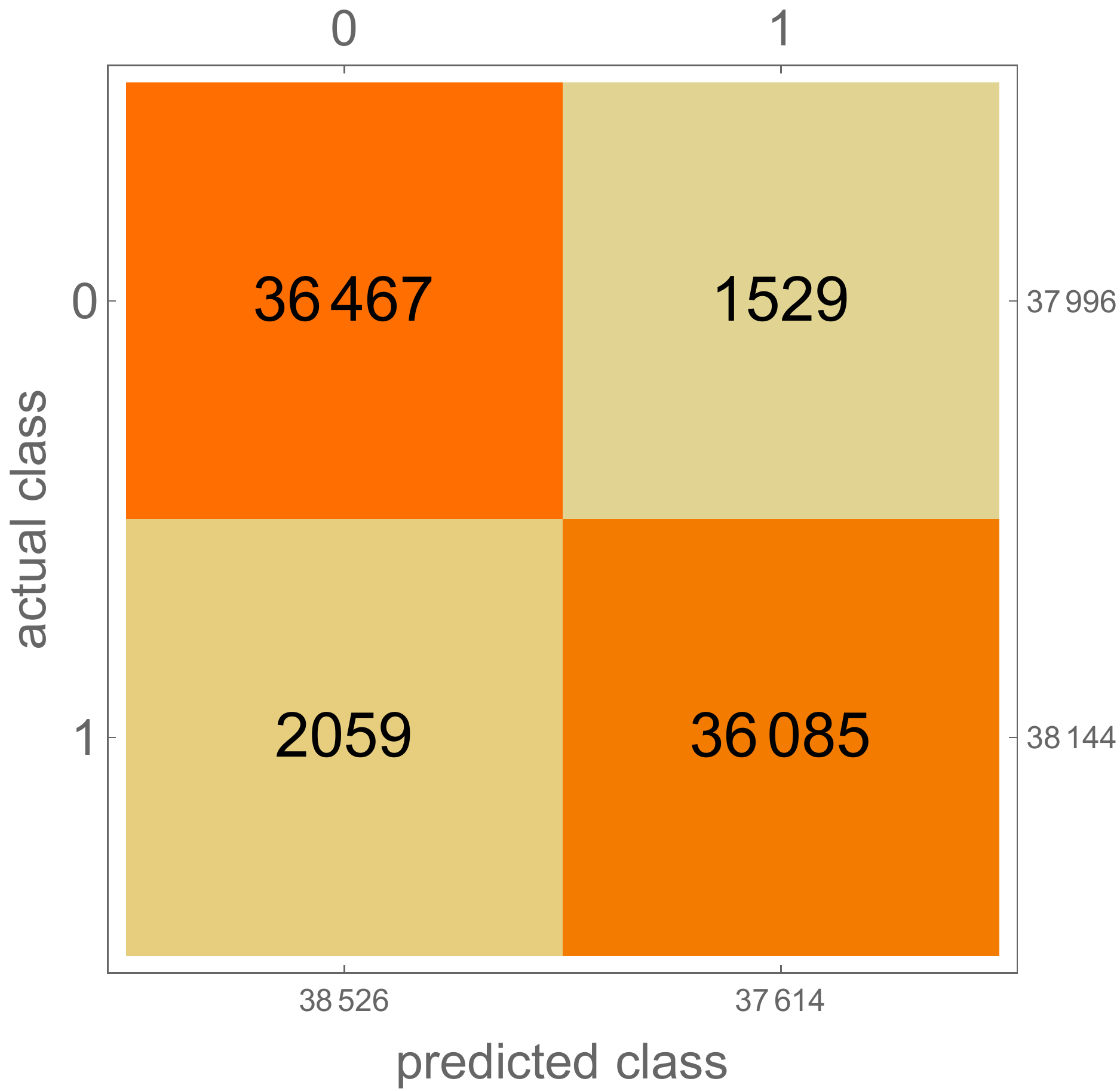}&\includegraphics[width=0.22 \textwidth]{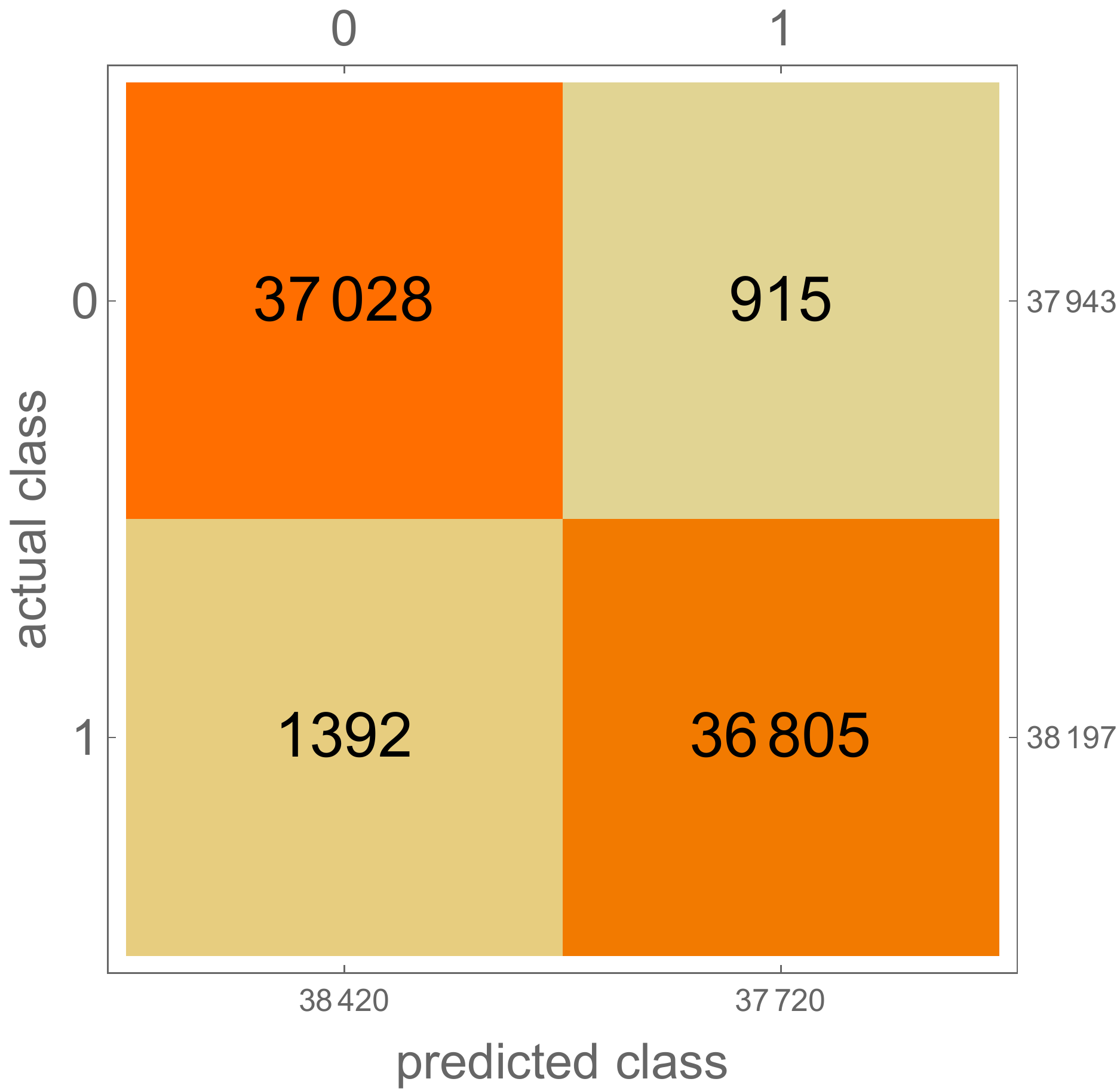} &\includegraphics[width=0.22 \textwidth]{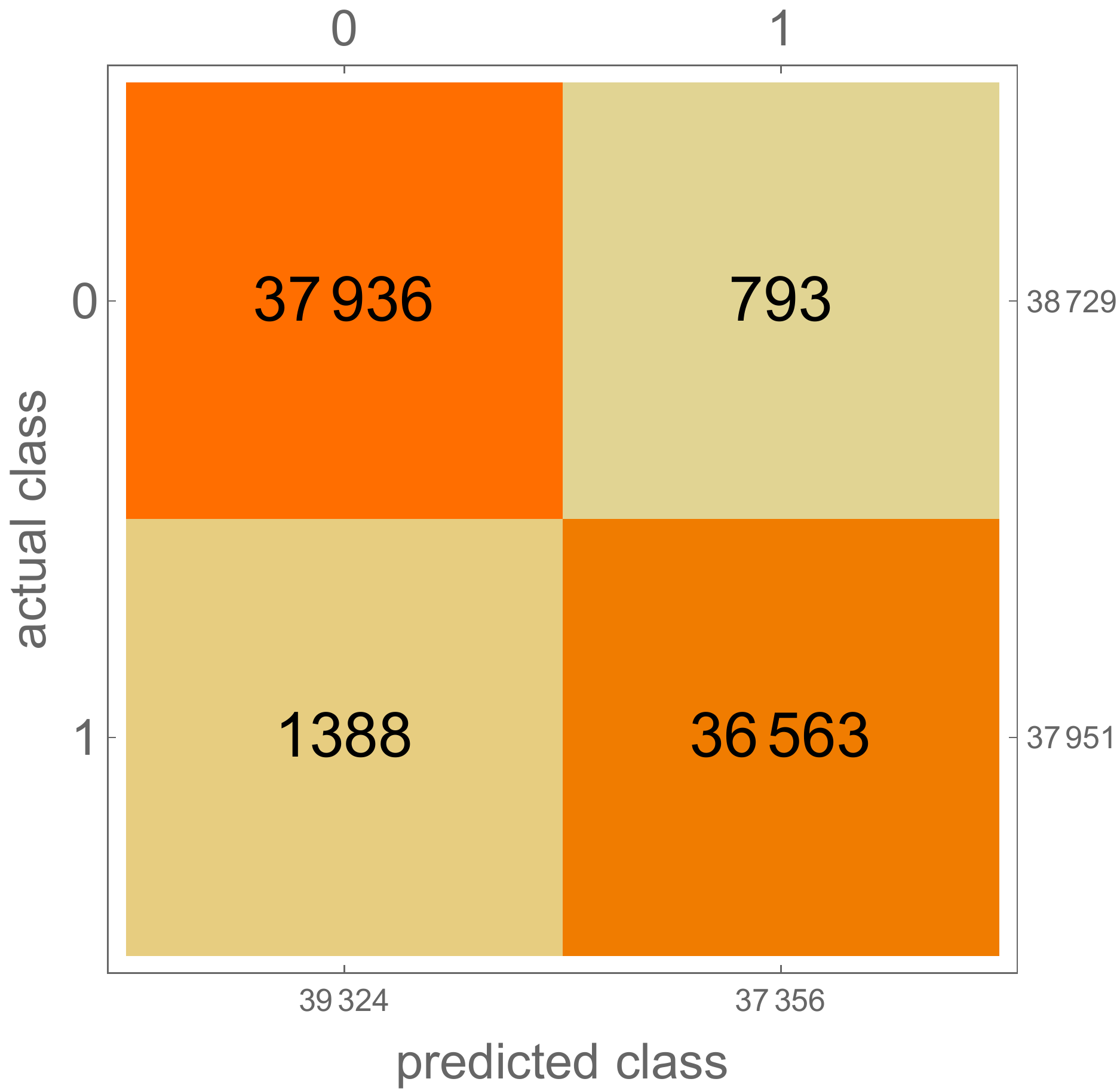}\\ \hline \raisebox{1.85 cm}{$13$}& 
\includegraphics[width=0.22 \textwidth]{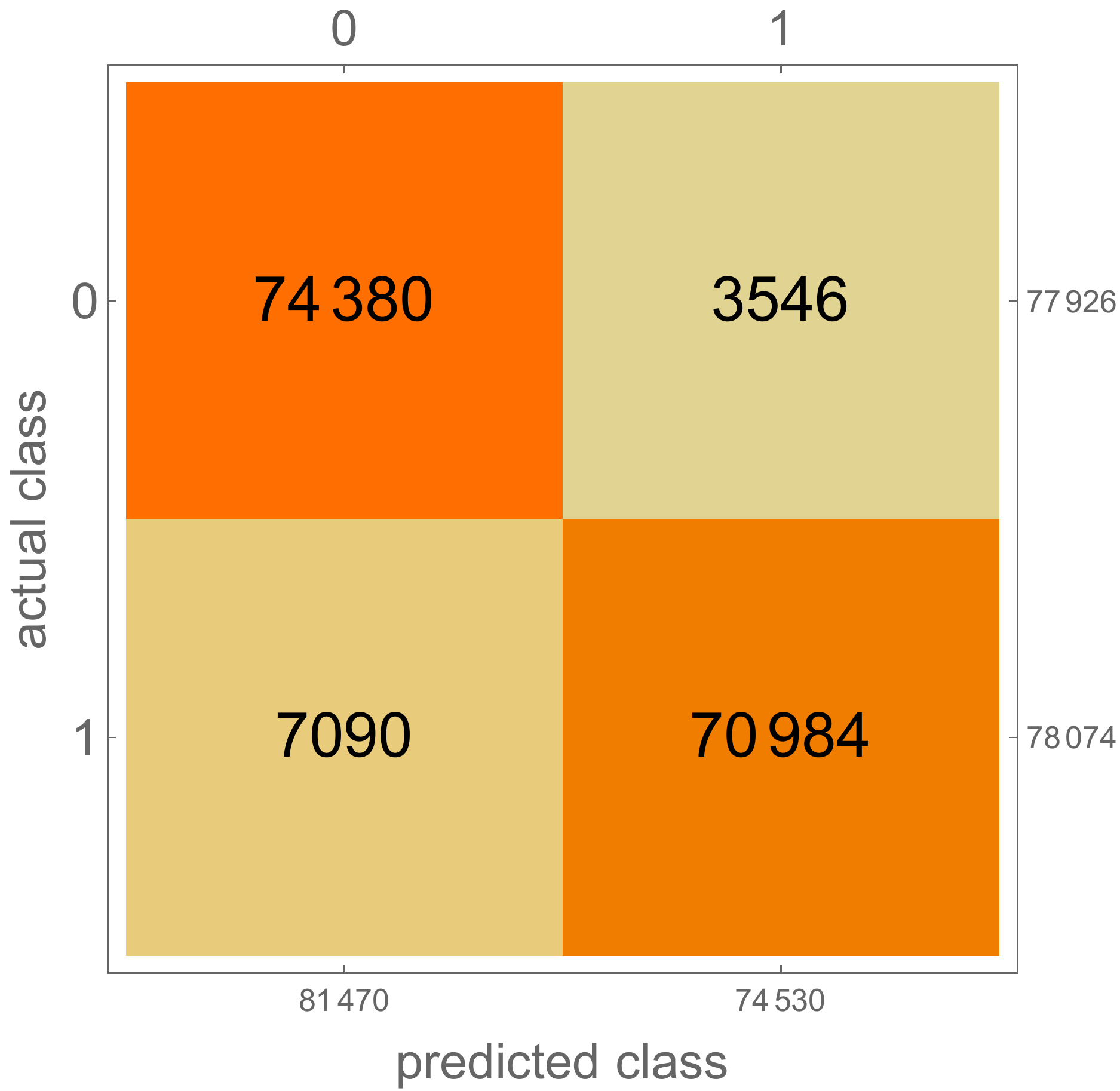} &\includegraphics[width=0.22 \textwidth]{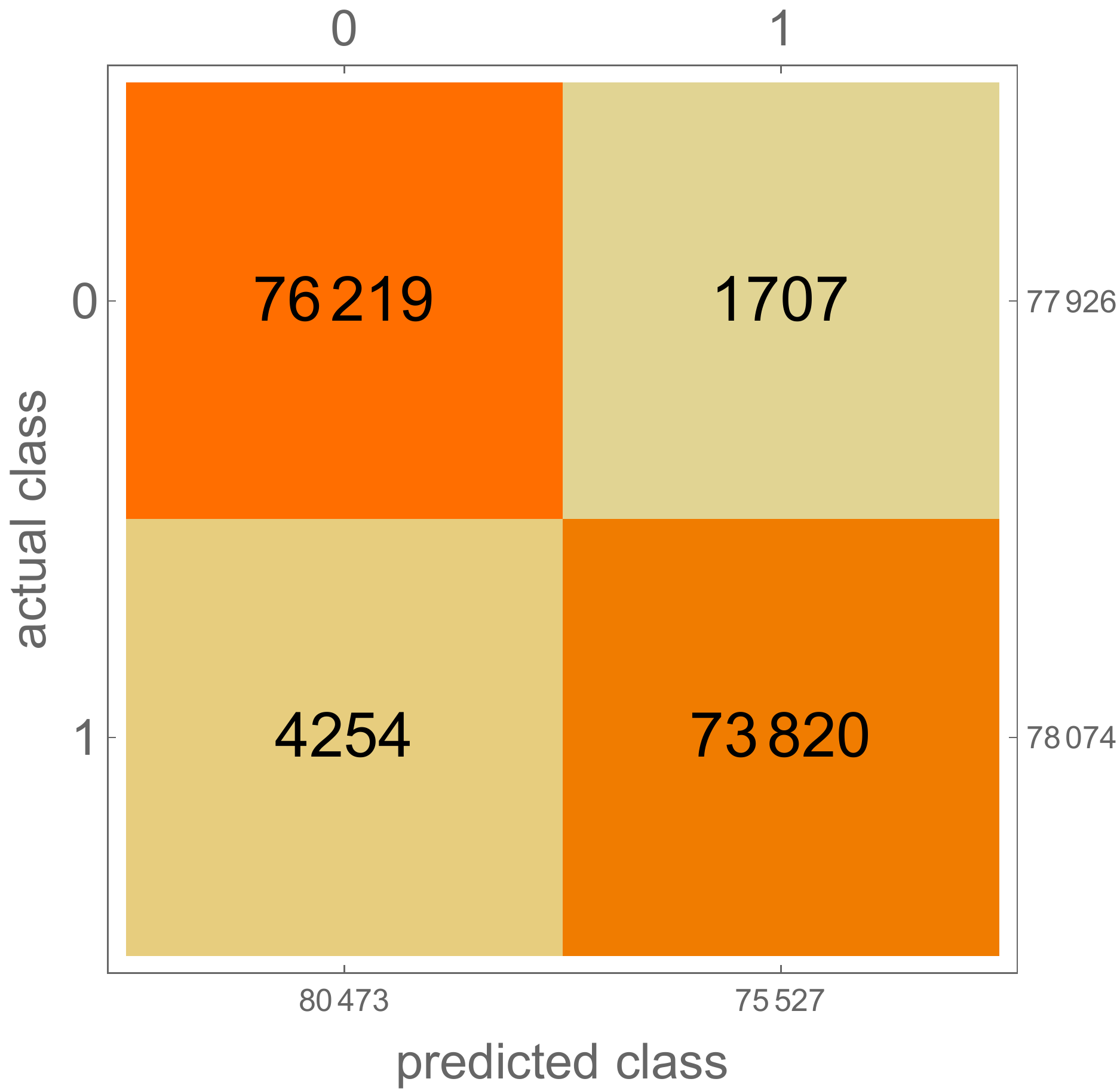}  &\includegraphics[width=0.22 \textwidth]{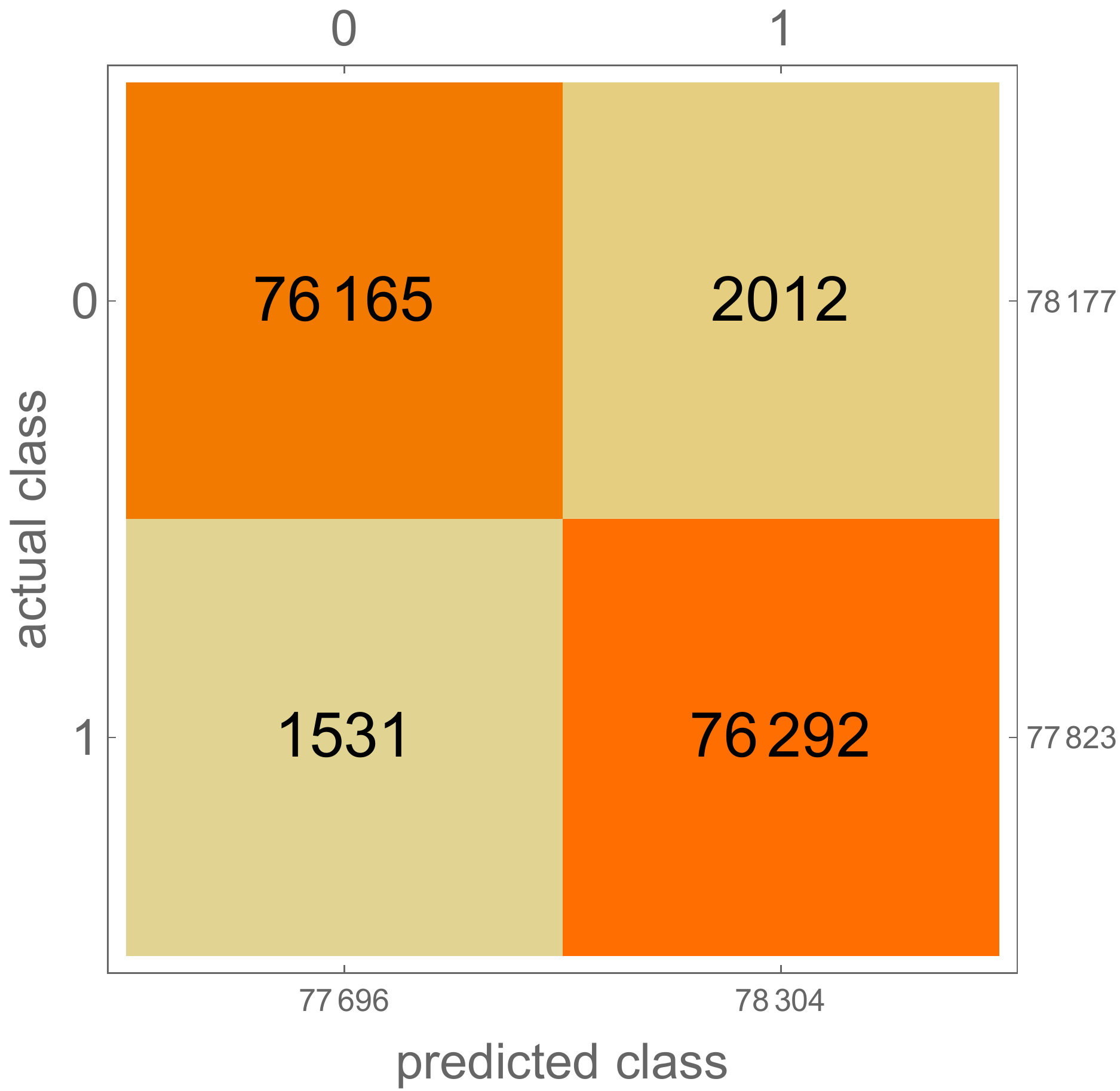}&\includegraphics[width=0.22 \textwidth]{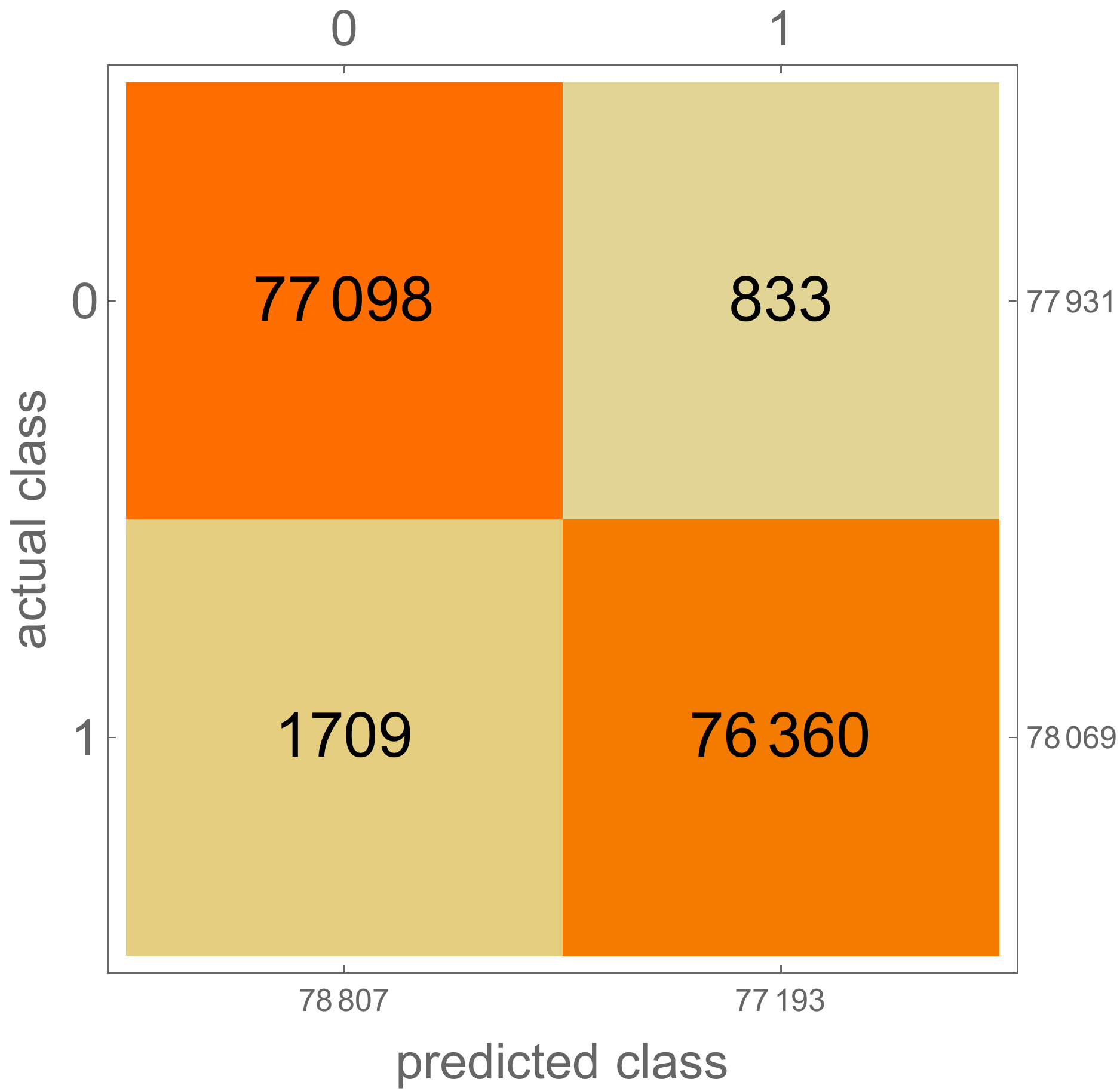}\\ \hline \raisebox{1.85 cm}{$14$} & 
\includegraphics[width=0.22 \textwidth]{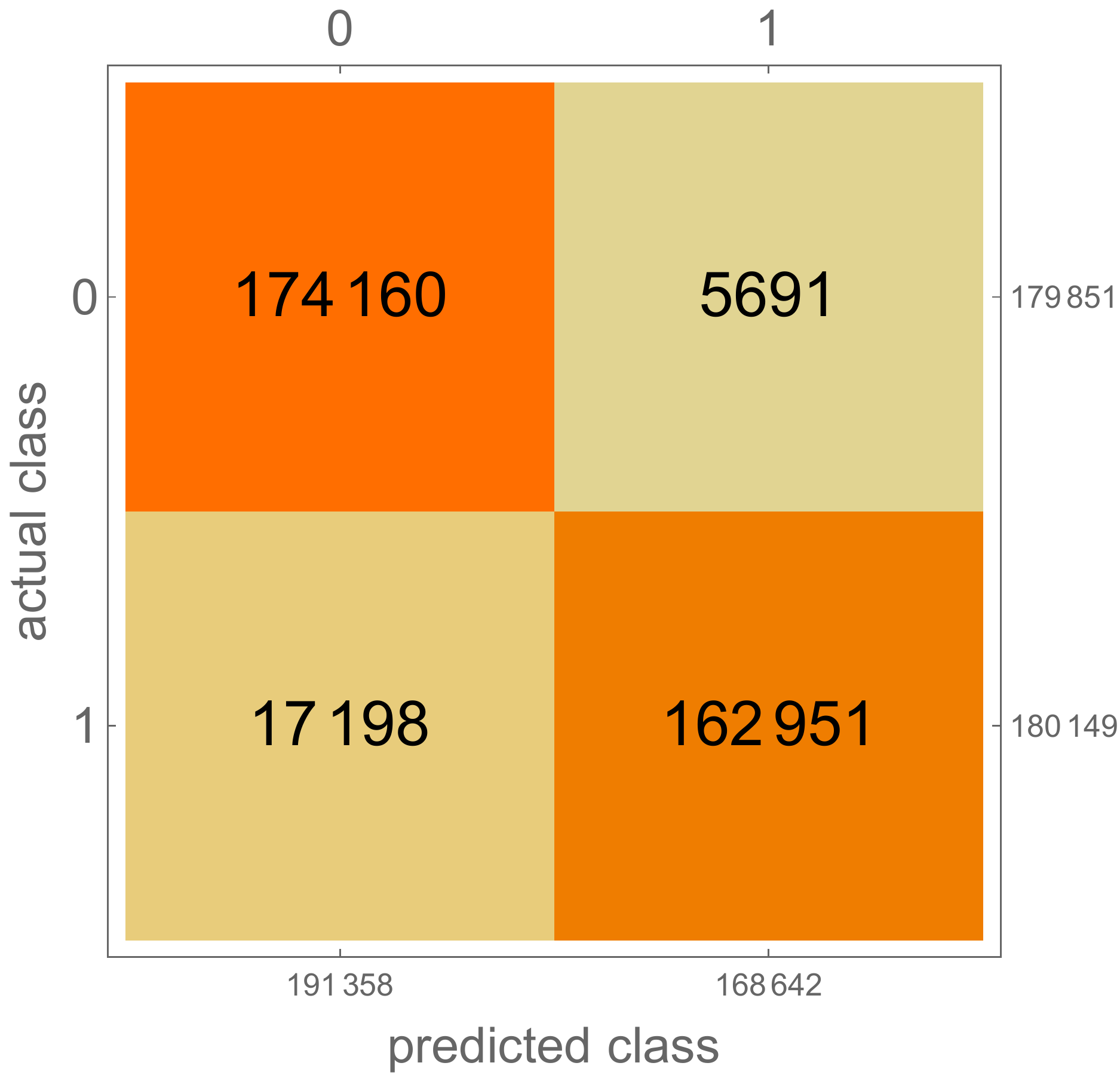} &\includegraphics[width=0.22 \textwidth]{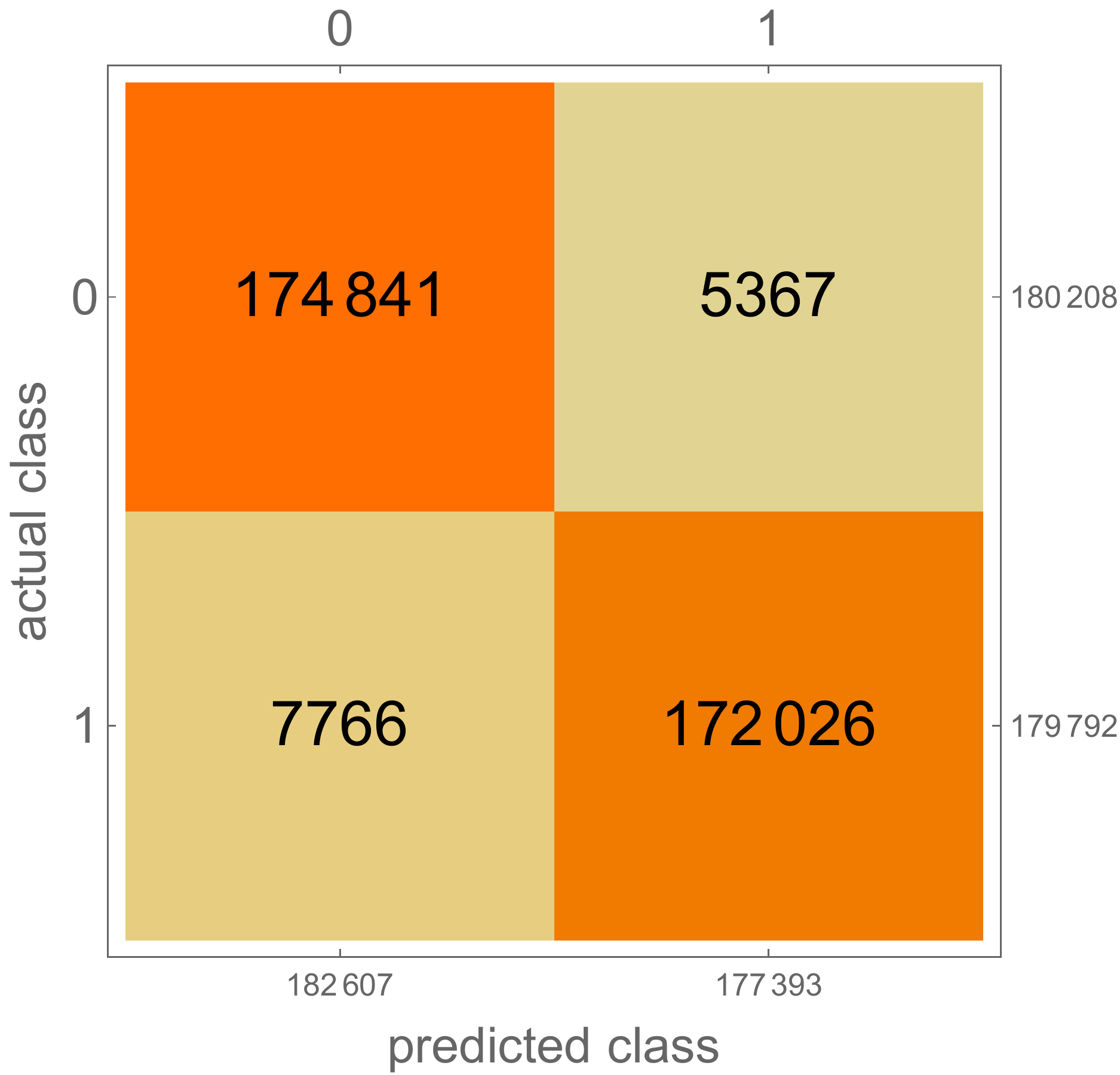}  &\includegraphics[width=0.22 \textwidth]{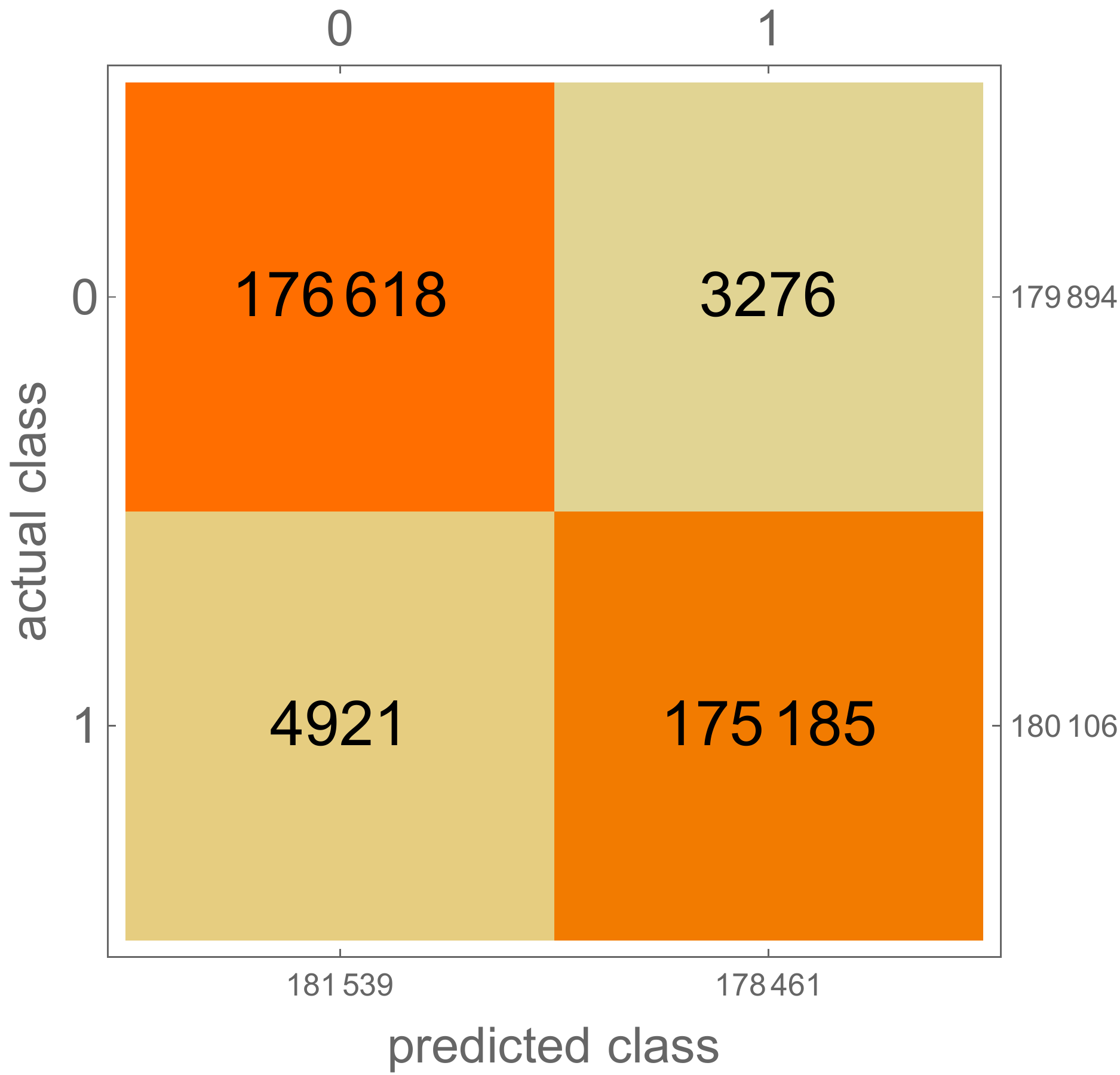}&\includegraphics[width=0.22 \textwidth]{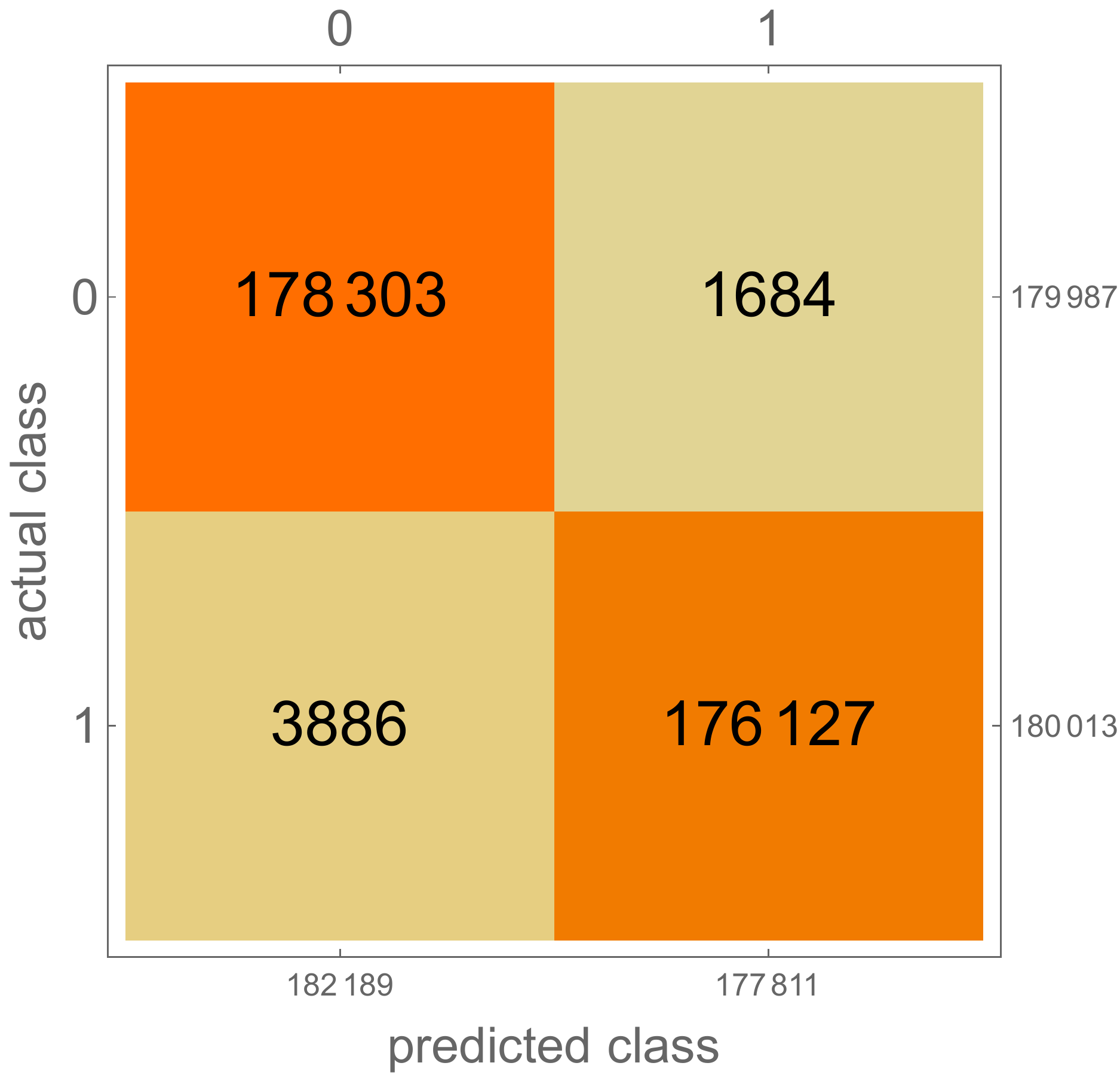}\\ \hline 
\end{tabular}
}
\end{center}
\medskip
\caption{\small \sf
The above table shows the confusion matrices of NearN, CNN$_2$, CNN$_3$ and LGBM classifiers in the experiments recorded in Table \ref{t:rank}. Here  class 0 represents the case $g_{\lambda, \mu}^\nu = 0$ and class 1 the case $g_{\lambda, \mu}^\nu \neq 0$.
}
\label{t:conf}
\end{table}

\section{Conclusions and Outlook}\label{s:outlook}

The experimental results in this paper show that an ML classifier can be trained to distinguish nonzero Kronecker coefficients from zero ones with high accuracy. This provides somewhat convincing evidence that ML can be used to recognize structures in Kronecker coefficients, adding another example to new opportunities of developing ML techniques for studying mathematical structures.

This paper focuses on the binary classification of $g_{\lambda, \mu}^\nu = 0$ or $g_{\lambda, \mu}^\nu >0$. One can perform several variations of this experiment. For example, multi-class classification of $g_{\lambda, \mu}^\nu=0,1, \dots , K$ or $g_{\lambda, \mu}^\nu > K$ can be considered for some $K \ge 1$. A different direction is to consider {\em induction} and {\em reduction} with respect to the inclusions $\mathfrak S_n \subset \mathfrak S_{N}$ for $n <N$. That is, we can train a machine with $\mathfrak S_{N}$ and apply it to $\mathfrak S_n$, and vice versa. Yet another direction is to investigate the datasets using dimensionality reduction techniques such as principal component analysis (PCA) or t-distributed stochastic neighborhood embedding (t-SNE) \cite{MH}.   

With all these experimental results at hand, the next step would be to perform attribute analysis to understand what ML classifiers actually recognize in the datasets. This may lead to some conjectures on how to combinatorially determine whether a Kronecker coefficient is zero or not, or how to combinatorially describe Kronecker coefficients. Hopefully, it may be that ML can even provide some hints about how to prove the conjectures.


\end{document}